\documentclass[12pt]{amsart}
\usepackage{amscd, amsfonts, amssymb}
\usepackage{fullpage}
\usepackage{latexsym}
\usepackage{verbatim}
\usepackage{enumerate}

\newcommand\bb{\mathfrak b}
\newcommand\cc{\mathfrak c}
\renewcommand\gg{\mathfrak g}
\newcommand\hh{\mathfrak h}
\newcommand\kk{\mathfrak k}
\newcommand\mm{\mathfrak m}
\newcommand\nn{\mathfrak n}
\newcommand\pp{\mathfrak p}
\newcommand\rr{\mathfrak r}
\newcommand\zz{\mathfrak z}

\newcommand{\inprod}[1]{\langle #1\rangle}

\newcommand\inverse{{^{-1}}}
\renewcommand{\check}{^{\vee}}

\newcommand\ra{\rightarrow}
\newcommand\lra{\longrightarrow}
\newcommand\sse{\subseteq}

\newcommand\UU{\mathcal U}
\newcommand\OO{\mathcal O}

\newcommand{\iso}{\cong}
\newcommand{\FF}{{\mathbb F}}
\newcommand{\NN}{{\mathbb N}}
\newcommand{\ZZ}{{\mathbb Z}}
\newcommand{\QQ}{{\mathbb Q}}
\newcommand{\RR}{{\mathbb R}}

\newcommand{\tuple}[1]{{\mathbf {#1}}}
\DeclareMathOperator{\ad}{ad}
\DeclareMathOperator{\Ad}{Ad}
\DeclareMathOperator{\Aut}{Aut}
\DeclareMathOperator{\Char}{char}
\DeclareMathOperator{\codim}{codim}
\DeclareMathOperator{\diag}{diag}
\DeclareMathOperator{\rk}{rk}
\DeclareMathOperator{\A}{A}
\DeclareMathOperator{\B}{B}
\DeclareMathOperator{\C}{C}
\DeclareMathOperator{\D}{D}
\DeclareMathOperator{\G}{G}
\DeclareMathOperator{\GL}{GL}
\DeclareMathOperator{\Gal}{Gal}
\DeclareMathOperator{\SL}{SL}
\DeclareMathOperator{\St}{St}
\DeclareMathOperator{\PGL}{PGL}
\DeclareMathOperator{\OR}{O}
\DeclareMathOperator{\Opp}{Opp}
\DeclareMathOperator{\SO}{SO}
\DeclareMathOperator{\SP}{Sp}
\DeclareMathOperator{\Hom}{Hom}
\DeclareMathOperator{\Lie}{Lie}
\DeclareMathOperator{\trdeg}{trdeg}
\DeclareMathOperator{\sat}{sat}
\DeclareMathOperator{\soc}{soc}
\DeclareMathOperator{\proj}{proj}
\DeclareMathOperator{\IM}{Im}
\DeclareMathOperator{\End}{End}
\DeclareMathOperator{\sgn}{sgn}
\DeclareMathOperator{\Mat}{Mat}
\DeclareMathOperator{\im}{im}

\numberwithin{equation}{section}

\newtheorem{thm}[equation]{Theorem}

\newtheorem{prop}[equation]{Proposition}

\theoremstyle{definition}

\newtheorem{exmp}[equation]{Example}

\theoremstyle{remark}
\newtheorem{rem}[equation]{Remark}
\theoremstyle{remark}

\newcommand{\ovl}{\overline}


\title[Complete reducibility and Steinberg endomorphisms]
{Complete reducibility and Steinberg endomorphisms}

\author[S. Herpel]{Sebastian Herpel}
\address
{Fakult\"at f\"ur Mathematik,
Ruhr-Universit\"at Bochum,
D-44780 Bochum, Germany}
\email{sebastian.herpel@rub.de}

\author[G. R\"ohrle]{Gerhard R\"ohrle}
\address
{Fakult\"at f\"ur Mathematik,
Ruhr-Universit\"at Bochum,
D-44780 Bochum, Germany}
\email{gerhard.roehrle@rub.de}

\begin{document}

\begin{abstract}
Let $G$ be a connected reductive algebraic group defined over  
an algebraically closed field of positive characteristic. 
We study a generalization of 
the notion of $G$-complete reducibility
in the context of Steinberg endomorphisms of $G$.
Our main theorem extends 
a special case of a rationality result in this setting.
\end{abstract}

\maketitle


\section{Introduction}
\label{sec:intro}

Let $p$ be a prime number and let $k=\ovl{\FF_p}$ be the algebraic closure of 
the field of $p$ elements.
Let $G$ be a connected reductive linear algebraic group defined over $k$ 
and let $H$ be a closed subgroup of $G$.
Let $\FF_p \subseteq k' \subseteq k$ be a field extension of $\FF_p$.
Following Serre \cite{serre2},
we say that a $k'$-defined subgroup $H$ of $G$ 
is  \emph{$G$-completely reducible over $k'$} provided that whenever 
$H$ is contained in a $k'$-defined parabolic subgroup $P$ of 
$G$, it is contained in a $k'$-defined Levi subgroup of $P$. 
If $k' = k$, then $H$ is $G$-completely reducible over $k'$ 
if and only if $H$ is $G$-completely reducible (or $G$-cr for short).
For an overview of this concept see for instance \cite{serre1} and 
\cite{serre2}.

The starting point for our discussion is the following special 
case of the rationality result \cite[Thm.\ 5.8]{BMR:2005}.
Let $q$ be a power of $p$ and let ${\FF}_q$ be the field of $q$ elements.

\begin{thm}
\label{thm:BMR5.8}
Suppose that both $G$ and $H$ are defined over $\FF_q$.
Then $H$ is $G$-completely reducible if and only if it is 
$G$-completely reducible over $\FF_q$.
\end{thm}

Let $\sigma: G \rightarrow G$ be a 
\emph{Steinberg endomorphism} of $G$, i.e.\ 
a surjective endomorphism of $G$ 
that fixes only finitely many points, see 
Steinberg \cite{Steinberg:1968End} for a detailed discussion
(for this terminology, see \cite[Def.\ 1.15.1b]{GLS:1998}). 
The set of all Steinberg endomorphisms of $G$ is a subset of all 
isogenies $G\rightarrow G$ (see \cite[7.1(a)]{Steinberg:1968End})
that encompasses in particular all (generalized) Frobenius endomorphisms, 
i.e.\ endomorphisms of $G$ some power of which are Frobenius endomorphisms 
corresponding to some $\FF_q$-rational structure on $G$. 

\begin{exmp}
Let $F_1,F_2$ be the Frobenius maps of $G=\SL_2$ given by raising coefficients
to the $p$th and $p^2$th powers, respectively.
Then the map 
$\sigma = F_1 \times F_2 : G \times G \rightarrow G \times G$
is a Steinberg morphism of $G \times G$ that is not a Frobenius morphism, cf.\
the remark following \cite[Thm.~2.1.11]{GLS:1998}.
\end{exmp}



If $G$ is almost simple, then $\sigma$ 
is a (generalized) Frobenius map 
(e.g.\ see \cite[Thm.\ 2.1.11]{GLS:1998}), 
and the possibilities for $\sigma$ are well known 
(\cite[\S 11]{Steinberg:1968End}, 
e.g.\ see \cite[Thm.\ 1.4]{Liebeck:1995}): 
$\sigma$ is conjugate to either $\sigma_q$, 
$\tau \sigma_q$, $\tau' \sigma_q$ or $\tau'$, 
where $\sigma_q$ is a standard Frobenius 
morphism, $\tau$ is an automorphism of 
algebraic groups coming from a graph 
automorphism of types $A_n$, $D_n$ or $E_6$, 
and $\tau'$ is a bijective endomorphism 
coming from a graph automorphism of type 
$B_2$ ($p=2$), $F_4$ ($p=2$) or $G_2$ ($p=3$). 

\begin{exmp}
If $G$ is not simple, then a generalized Frobenius map may fail to factor 
into a field and a graph automorphism as stated above. 
For example, let $p=2$ and 
let $H_1, H_2$ be simple, simply connected groups of type $B_n$
and $C_n$ ($n \ge 3$), respectively.
Then there are special isogenies
$\phi_1 : H_1 \rightarrow H_2$ and $\phi_2: H_2 \rightarrow H_1$ whose
composites $\phi_1 \circ \phi_2$ and $\phi_2 \circ \phi_1$ are standard
Frobenius maps with respect to $p$ on $H_2$, respectively $H_1$,
see \cite[p 5 of Exp.\ 24]{Chevalley:2005}.
Let $G = H_1 \times H_2$ and define $\sigma: G \rightarrow G$ by
$\sigma(h_1,h_2) = (\phi_2(h_2),\phi_1(h_1))$.
Then $\sigma$ 
is an example of such a more complicated generalized Frobenius map.
\end{exmp}

We now give an extension of Serre's notion of 
$G$-complete reducibility in this setting of 
Steinberg endomorphisms:
Let $\sigma$ be a Steinberg endomorphism of $G$ 
and let $H$ be a 
subgroup of $G$. We say that $H$ is 
\emph{$\sigma$-completely reducible} (or 
$\sigma$-cr for short), provided that whenever 
$H$ lies in a $\sigma$-stable parabolic subgroup 
$P$ of $G$, it lies in a $\sigma$-stable Levi subgroup of $P$.
This notion is motivated as follows: 
If $\sigma_q$ is a standard Frobenius morphism of $G$, 
then a subgroup $H$ of $G$ is defined over $\FF_q$ if and only 
if it is $\sigma_q$-stable and if so, 
$H$ is $G$-completely reducible 
over $\FF_q$ if and only if it is $\sigma_q$-completely reducible.
In view of this new notion, 
the goal of this note is the following generalization of 
Theorem \ref{thm:BMR5.8} to arbitrary
Steinberg endomorphisms of $G$ (the special case of Theorem \ref{thm:main} 
when $\sigma = \sigma_q$ gives Theorem \ref{thm:BMR5.8}). 

\begin{thm}
\label{thm:main}
Let 
$\sigma$ be a Steinberg endomorphism of $G$. 
Let $H$ be a $\sigma$-stable subgroup 
of $G$. Then $H$ is $G$-completely reducible 
if and only if $H$ is $\sigma$-completely reducible. 
\end{thm}

Theorem \ref{thm:main} follows from 
Theorems \ref{thm:sigma->G-cr} and \ref{thm:G-cr->sigma}
proved in the next section. 

\begin{exmp}
Theorem \ref{thm:main} is false without the 
$\sigma$-stability condition on $H$.
For instance, a maximal torus $T$ of $G$ is always $G$-cr, 
cf.\ \cite[Lem.\ 2.6]{BMR:2005}. 
But it may happen that $T$ is contained in a $\sigma$-stable Borel subgroup of $G$,
without being itself $\sigma$-stable. Then $T$ clearly fails to be $\sigma$-cr.
In the other direction, $G$ may contain a maximal parabolic subgroup 
$P$ of $G$ that is not $\sigma$-stable. 
The only $\sigma$-stable parabolic subgroup of $G$ containing $P$ is $G$ itself.
Then $P$ is $\sigma$-cr for trivial reasons, whereas a proper
parabolic subgroup of $G$ is not $G$-cr.
\end{exmp}

\begin{rem}
Even if $H$ is not $\sigma$-stable, Theorem \ref{thm:main}  
gives some information about the notion of $\sigma$-complete reducibility, 
as follows. Let $\overline{H}^\sigma$ be the algebraic 
subgroup of $G$ generated by all translates $\sigma^i H$, $i \geq 0$. Then
$\overline{H}^\sigma$ is $\sigma$-stable and contained in the same 
$\sigma$-stable subgroups of $G$ as $H$. 
In particular, $H$ is $\sigma$-cr if and
only if $\overline{H}^\sigma$ is $\sigma$-cr. Thus, by Theorem \ref{thm:main}, 
this is equivalent to $\overline{H}^\sigma$ being $G$-cr. 
\end{rem}


\section{Proof of  Theorem \ref{thm:main}}
\label{sect2}

In addition to the notation already fixed in the Introduction, 
$\sigma : G \to G$
is always a Steinberg endomorphism
of $G$ and from now on the subgroup $H$ of $G$ is assumed to be 
$\sigma$-stable.
We begin with a generalization of (a special case of) 
\cite[Prop.\ 2.2 and Rem.\ 2.4]{LMS:2005}.

\begin{prop} 
\label{prop:dichotomy}
If $H$ is not $G$-completely reducible, then there 
exists a proper $\sigma$-stable parabolic subgroup 
of $G$ containing $H$.
\end{prop}

\begin{proof}
First we assume that $G$ is almost simple. 
We want to reduce to the case where $H$ is a finite, 
$\sigma$-stable subgroup of $G$, and then apply 
\cite[Prop.\ 2.2 and Rem.\ 2.4]{LMS:2005}. 
Since $G$ is almost simple, we can assume that $\sigma^m = \sigma_q$ 
is a standard Frobenius map for some positive integer $m$. 
We choose a closed embedding $G \rightarrow \GL_n(k)$ so that 
$\sigma_q$ is the restriction of the standard Frobenius map 
of $\GL_n(k)$ that raises coefficients to the $q$th power 
(see \cite[Prop.\ 4.1.11]{Geck:2003}). 
For $r \in \ZZ, r \geq 1$, 
let $\tilde{H}(r) = H \cap \GL_n({\FF_{q^{r!}}})$. 
Then we can write $H$ as the directed union of finite 
subgroups $H = \bigcup_{r\geq 1} \tilde{H}(r)$. 
Note that the union is indeed directed, that is
\begin{align} 
\label{eq:inclusion}
\tilde{H}(r) \subseteq \tilde{H}(r+1) ~ \forall r\geq 1.
\end{align}
We wish to construct a similar, but $\sigma$-stable filtration of $H$. 
For this purpose we set $H(r) = \bigcap_{l=0}^{m-1} \sigma^l \tilde{H}(r)$. 
Then each $H(r)$ is a finite, $\sigma$-stable subgroup of $H$ 
(for the $\sigma$-stability, we use that each $\tilde{H}(r)$ 
is stable under $\sigma^m = \sigma_q$). 
Moreover, we claim that $H$ is the directed union 
$H = \bigcup_{r\geq 1} H(r)$. 
Indeed, if $h \in H$, then the identities 
$H=\sigma H$ and $H = \bigcup_{r \geq 1} \tilde{H}(r)$ 
imply that for each $l=0,\dots,m-1$ we can find some 
$r_l$ such that $h \in \sigma^l \tilde{H}(r_l)$. 
But then \eqref{eq:inclusion} implies that 
$h \in H(r)$ for $r \geq \max\{r_0,\dots,r_{m-1}\}$. 
It follows from the argument in the proof of \cite[Lem.\ 2.10]{BMR:2005}
that there is an integer $r'$ so that $H(r')$ has the following property: 
$H$ is contained in a parabolic subgroup 
$P$ of $G$ (respectively a Levi subgroup $L$ of $G$) 
if and only if $H(r')$ is contained in $P$ (respectively in $L$). 
Therefore, if $H$ is not $G$-cr, then neither is $H(r')$, 
and we can apply \cite[Prop.\ 2.2 and Rem.\ 2.4]{LMS:2005}  
to obtain a proper $\sigma$-stable parabolic subgroup $P$ of 
$G$ that contains $H(r')$.  But then $P$ also contains $H$.

Next we drop the simplicity assumption on $G$. 
Then we can use the almost simple components of $G$ 
to reduce to the almost simple case: 
Let $\pi:G' := Z(G)^\circ \times G_1 \times \dots \times G_r \rightarrow G$ 
be the product map, where $G_1, \dots, G_r$ 
are the almost simple components of the semisimple group 
$[G,G]$ and let $\pi_i:G' \rightarrow G_i$ be the projection ($1\leq i\leq r$). 
Then $\pi$ is an isogeny. Let $H' = \pi^{-1}(H)$. 
Using \cite[Lem.\ 2.12]{BMR:2005} and the fact that 
$Z(G)^\circ$ is a torus, 
we find that there is some 
index $i$ such that $H_i := \pi_i(H') \subseteq G_i$ is not $G_i$-cr. 
We can assume that $i=1$. 
We are now in the situation of the first part of the proof (for $H_1 \subseteq G_1$), 
except that we have yet to specify a Steinberg endomorphism 
of $G_1$ that stabilizes $H_1$. 
Since $\sigma$ stabilizes $[G,G]$ 
and maps components to components (\cite[Exp.\ 18, Prop.\ 2]{Chevalley:2005}), 
we can assume that $\sigma$ permutes $G_1,\dots,G_s$ 
cyclically for some $s \leq r$. 
Moreover, $\sigma$ stabilizes $Z(G)^\circ = R(G)$ 
(because $\sigma$ is an isogeny). 
Using the restrictions $\sigma|_{Z(G)^\circ}$ and $\sigma|_{[G,G]}$, 
we can define a Steinberg endomorphism $\sigma':G' \rightarrow G'$ 
of $G'$ such that $\pi \circ \sigma' = \sigma \circ \pi$.
We denote by $H''$ the image (under the projection) 
of $H'$ in $G'' := G_1 \times \dots \times G_s$. 
Now let $\tau = \sigma^s|_{G_1}:G_1 \rightarrow G_1$ 
denote the generalized Frobenius map on $G_1$ 
induced by $\sigma$ (\cite[Thms.\ 2.1.2(g) and 2.1.11]{GLS:1998}). 
Then $H_1$ is $\tau$-stable, since $H$ is $\sigma^s$-stable. 
We apply the first part of the proof to $H_1 \subseteq G_1$ 
to obtain a proper $\tau$-stable parabolic subgroup $P_1$ of $G_1$ 
containing $H_1$. 
Then $P'':= P_1 \times \sigma P_1 \times \dots \times \sigma^{s-1}P_1 \subseteq G''$ 
is a proper $\sigma'|_{G''}$-stable parabolic subgroup of $G''$ 
(\cite[Cor.\ 6.2.8]{Springer:1998}). 
The bijectivity of $\sigma^s|_{H_i}:H_i \rightarrow H_i$ 
for $1 \leq i \leq s$ implies that $H_i = \sigma^{i-1}H_1$ 
for $1 \leq i \leq s$. We get that $P''$ contains $H''$, 
since we have  
$H'' \subseteq H_1 \times H_2 \times\dots \times H_s$ 
and $H_1\subseteq P_1$. 
Consequently, 
$P' = Z(G)^\circ \times P'' \times G_{s+1} \times \dots \times G_r$ 
is a proper $\sigma'$-stable parabolic subgroup of $G'$ containing $H'$. 
Finally, $P = \pi(P')$ is a proper 
$\sigma$-stable parabolic subgroup of $G$ containing $H$, as desired.
\end{proof}

\begin{rem}
In \cite[Prop.\ 2.2 and Rem.\ 2.4]{LMS:2005}, 
Liebeck, Martin and Shalev prove the following: 
Let $G$ be an almost simple algebraic group over $k$ as above. 
Let $\Aut^\#(G)$ denote the group generated by inner automorphisms of $G$, 
together with $p^i$-power field morphisms ($i \geq 1$), 
and graph automorphisms (which may include the bijective 
endomorphisms coming from a graph automorphism of type 
$B_2$ ($p=2$), $F_4$ ($p=2$) or $G_2$ ($p=3$)). 
(Note that $\Aut^\#(G)$ is an extension of the group 
$\Aut^+(G)$ from \cite{LMS:2005}.)
Let $S$ be a subgroup of $\Aut^\#(G)$ and suppose 
that $H \subseteq G$ is a finite, $S$-stable subgroup 
that is not $G$-cr. Then $H$ is contained in a proper 
$S$-invariant parabolic subgroup of $G$ 
(note that the notion of strongly reductive subgroups in $G$ 
is equivalent to the notion of $G$-completely reducible subgroups, 
cf.\ \cite[Thm.\ 3.1]{BMR:2005}). 
If we take $S$ to be generated by a (generalized) 
Frobenius endomorphism $\sigma$ of $G$, 
then we get the assertion of 
Proposition \ref{prop:dichotomy} for $G$ almost simple and $H$ finite.  
\end{rem}

\begin{thm} 
\label{thm:sigma->G-cr}
If $H$ is $\sigma$-completely reducible, 
then it is $G$-completely reducible.
\end{thm}

\begin{proof}
If $H$ is not contained in any 
proper $\sigma$-stable parabolic subgroup of $G$, 
then it is $G$-cr according to Proposition \ref{prop:dichotomy}. 
So we can assume that there is a proper $\sigma$-stable 
parabolic subgroup $P$ of $G$ containing $H$. 
We choose $P$ minimal with these properties. 
Since $H$ is $\sigma$-cr, it is contained 
in a $\sigma$-stable Levi subgroup $L$ of $P$. 
Suppose there is a proper $\sigma$-stable 
parabolic subgroup $P_L$ of $L$ containing $H$. 
Then $P'=P_L R_u(P) \subsetneq P$ is another 
parabolic subgroup of $G$ 
(see \cite[Prop.\ 4.4(c)]{BT:1965}) containing $H$, 
and $P'$ is $\sigma$-stable ($\sigma$ stabilizes $R_u(P)$ as any isogeny does). 
But this contradicts our choice of $P$. 
So we can use Proposition \ref{prop:dichotomy} 
again to deduce that $H$ is $L$-cr, which in 
turn implies that $H$ is $G$-cr (\cite[Cor.\ 3.22]{BMR:2005}).     
\end{proof}

For the converse of Theorem \ref{thm:sigma->G-cr} 
we argue as in the last part of the proof of \cite[Thm.\ 9]{LS:1998}. 
But first we recall a parametrization of 
the parabolic and Levi subgroups of $G$ in terms of cocharacters of $G$,
e.g.\ see  \cite[Lem.\ 2.4]{BMR:2005}: 
Given a parabolic subgroup $P$ of $G$ and any Levi subgroup $L$ of $P$, 
there exists some cocharacter $\lambda$ of $G$ such that $P$ and $L$ 
are of the form 
$P = P_\lambda = \{g \in G \mid \lim_{t\rightarrow 0} \lambda(t)g\lambda(t)^{-1} \text{ exists}\}$ 
and $L = L_\lambda = C_G(\lambda(k^*))$, respectively. 
The unipotent radical of $P_\lambda$ is then given by 
$R_u(P_\lambda) = \{g \in G \mid \lim_{t\rightarrow 0} \lambda(t)g\lambda(t)^{-1} = 1\}$.

\begin{thm}
\label{thm:G-cr->sigma}
If $H$ is $G$-completely reducible,
then it is $\sigma$-completely reducible.
\end{thm}

\begin{proof}
Suppose that $P$ is a $\sigma$-stable 
parabolic subgroup of $G$ containing $H$. 
Since $H$ is $G$-cr, there is some Levi 
subgroup $L$ of $P$ that contains $H$. 
Let $U=R_u(P)$. 
Then $\Lambda = \{uLu^{-1} \mid u \in U, H \subseteq uLu^{-1}\}$ 
is the set of all Levi subgroups of $P$ that contain $H$. 
Clearly, $\Lambda$ is $\sigma$-stable, since $H$ and $P$ are. 
We need to prove that $\Lambda$ contains an element fixed by $\sigma$. 

If $uLu^{-1}$ is in $\Lambda$, 
then $u^{-1}Hu \subseteq L \cap UH = H$, 
so that $u$ normalizes $H$. 
In fact, $u$ centralizes $H$, 
since $[N_U(H),H]\subseteq H\cap U = \{1\}$. 
So the group $C = C_U(H)$ acts transitively on $\Lambda$. 
We claim that $C$ is connected. In order to prove this, 
write $P = P_\lambda$, $L=L_\lambda$ and $U = R_u(P_\lambda)$ 
for some suitable cocharacter $\lambda$ of $G$. 
The torus $\lambda(k^*)$ normalizes $C_G(H)$ 
(because $H$ is contained in $L$) and $U$, 
hence it normalizes $C$. 
Whence, for any fixed $c\in C$,   
the map $\phi_c : k^*\rightarrow C$, given by 
$t \mapsto \lambda(t)c\lambda(t)^{-1}$, is well-defined. 
Moreover, $C \subseteq U$ implies that $\phi_c$ extends 
to a morphism $\hat \phi_c : k\rightarrow C$ that maps $0$ 
to $1$ and $1$ to $c$. Since the image of $\hat \phi_c$ 
is connected, we get $c \in C^\circ$. 
It follows that $C = C^\circ$.
But now we can apply the 
Lang-Steinberg theorem (see \cite[Thm.\ 10.1]{Steinberg:1968End}) 
to conclude that $\Lambda$ contains an element fixed by $\sigma$.
\end{proof}

\begin{rem}
We conclude by outlining a short alternative approach to  
Proposition \ref{prop:dichotomy}; the latter was crucial in the proof of 
Theorem \ref{thm:sigma->G-cr}.
This variant utilizes 
the so called \emph{Centre Conjecture} for spherical buildings due to J.~Tits 
from the 1950s.
This deep conjecture has recently been established by work of
Leeb and Ramos-Cuevas, e.g.\ see \cite[\S 2]{BMR:2010} 
and the references therein for further details.
This conjecture states that in the building $\Delta = \Delta(G)$ of $G$ 
any convex contractible subcomplex $\Sigma$ has a simplex which is fixed under 
any building 
automorphism of $\Delta$ which stabilizes $\Sigma$ as a subcomplex. Such 
a fixed simplex is often referred to as a \emph{centre} giving this conjecture 
its name. 
Here is a sketch of a building theoretic alternative to 
the proof of Proposition \ref{prop:dichotomy}:
Let $H$ be a $\sigma$-stable subgroup of $G$
which is not $G$-cr.
Consider the subcomplex $\Delta^H$ of $H$-fixed points of the building 
$\Delta$, i.e., $\Delta^H$ corresponds to the set of all parabolic subgroups
of $G$ that contain $H$. Note that $\Delta^H$ 
is always convex (\cite[Prop.\ 3.1]{serre2}) 
and since $H$ is not $G$-cr, $\Delta^H$ is also 
contractible (\cite[Thm.\ 2]{serre1.5}).
The Steinberg morphism $\sigma$ of $G$ affords 
a building automorphism of $\Delta$, also denoted by $\sigma$.
Since $H$ is $\sigma$-stable, so is $\Delta^H$.
Now since $\Delta^H$ is convex and contractible, the Centre Conjecture
asserts the existence of a centre of $\Delta^H$ with respect to the action of $\sigma$ 
which corresponds to a proper parabolic subgroup of $G$ 
which is $\sigma$-stable and contains $H$.
This is precisely the conclusion of Proposition \ref{prop:dichotomy}.
\end{rem}



{\bf Acknowledgments}:
The authors acknowledge the financial support of
the DFG-priority program SPP 1388 ``Representation Theory''.
We are grateful to Olivier Brunat for helpful discussions on the material
of this note.


\end{document}